\documentclass{llncs}
\usepackage[latin1]{inputenc} 
\usepackage[T1]{fontenc}      
\usepackage[dvips]{graphicx}
\usepackage{epsf}
\usepackage{latexsym}
\usepackage{epsfig}
\usepackage{amsmath}
\usepackage{amsfonts}
\usepackage{amssymb}
\usepackage{algorithm}
\usepackage{alltt}

\def\R{\mbox{\ensuremath{I \! \! R}}}

\def\Xdo{\mbox{\bf do }}
\def\Xwhile{\mbox{\bf while }}

\def\Xif{\mbox{\bf if }}
\def\Xthen{\mbox{\bf then }}
\def\Xelse{\mbox{\bf else }}

\def\Xendwhile{\mbox{\bf endwhile }}

\def\Xreturn{\mbox{\bf return }}

\def\Xaffect{{\mbox{\ensuremath{\leftarrow}}}}

\newcommand{\lb}[1]{\ensuremath{\underline{#1}}}
\newcommand{\ub}[1]{\ensuremath{\overline{#1}}}

\newcommand{\vi}[1]{\ensuremath{\text{\bf #1}}}

\newcommand{\calx}[1]{\ensuremath{{\mathcal #1}}}
\newcommand{\Pb}{\ensuremath{{\mathcal P}}}

\newcommand{\url}[1]{{\it #1}}

\newcommand{\vnull}{{0}}

\newcommand{\vg}{{g}}

\newcommand{\vl}{{l}}

\newcommand{\vx}{{x}}

\begin{document}

\title{Revisiting the upper bounding process in a safe Branch and
  Bound algorithm\protect\footnote{An
    extented version of this paper is available  at:\\
  http://www.i3s.unice.fr/\%7Emh/RR/2008/RR-08.11-A.GOLDSZTEJN.pdf}} 
\author{Alexandre Goldsztejn\inst{1} \and Yahia Lebbah\inst{2,3} \and Claude Michel\inst{3} \and Michel Rueher\inst{3}}

\institute{CNRS / Université de Nantes 2, rue de la Houssinière, 44322 Nantes, France
\email{alexandre.goldsztejn@univ-nantes.fr}\\
\and
Universit\'e d'Oran Es-Senia B.P. 1524 EL-M'Naouar, 31000 Oran, Algeria
\email{ylebbah@gmail.com}\\
\and
Universit{\'e} de Nice-Sophia Antipolis, I3S-CNRS,  06903 Sophia Antipolis, France
\email{\{cpjm, rueher\}@polytech.unice.fr}}  

\maketitle

\begin{abstract}
Finding feasible points for which the proof succeeds is a critical
issue in safe   Branch and  Bound algorithms   which handle continuous problems. 
In this paper, we introduce a new strategy to compute 
very accurate approximations of feasible points. 
This strategy takes advantage of the Newton 
method for under-constrained systems of equations and inequalities. More precisely, it
exploits the optimal solution of a linear relaxation of the problem to compute 
efficiently a promising upper bound.
First experiments on the Coconuts benchmarks demonstrate that this approach is very effective. 
\end{abstract}


\subsection*{Introduction}
Optimization problems are a challenge for CP in finite domains; they
are also a big challenge for CP on continuous domains.
The point is that CP solvers are much  slower than classical
(non-safe) mathematical methods on nonlinear constraint problems as
soon as we consider optimization problems. The techniques introduced
in this paper try to boost constraints techniques on these problems and
thus, to reduce the gap between efficient but unsafe systems like
BARON\footnote{See http://www.andrew.cmu.edu/user/ns1b/baron/baron.html},
and the slow but safe constraint based approaches.  
We consider here the global optimization problem {\Pb} to minimize an
objective function under nonlinear equalities and inequalities,

\begin{equation}\label{equ-nonlinprog}
\begin{array}{ll}
$minimize$ & f(x)\\
$subject to$ & g_i(x) = 0,\,\, i \in \{1, .., k\} \\
             & h_j(x) \leq 0,\,\, j \in \{1, .., m\}
\end{array}
\end{equation}

with $x \in \vi{x}$, $f:\R^n\rightarrow \R$, $g_i:\R^n\rightarrow \:\R$ and $h_j:\R^n\rightarrow \:\R$;
Functions $f$ , $g_i$ and $h_j$ are nonlinear and continuously differentiable on some
vector $\vi{x}$  of intervals of $\R$.
For convenience, in the sequel, $g(x)$ (resp. $h(x)$) will denote the vector of
$g_i(x)$ (resp. $h_j(x)$) functions. 

The difficulties in such global optimization problems come mainly from the fact that
many local minimizers may exist but only few  of them are global minimizers \cite{Neumaier-2004}.
Moreover, the feasible region may be disconnected. 
Thus, finding feasible points is a critical issue in safe Branch and 
Bound algorithms for continuous global optimization. 
Standard strategies use local search techniques to provide a reasonable approximation of an upper bound
and try to prove that a feasible solution actually exists within the
box around the guessed global optimum. 
Practically, finding a guessed point for which the proof succeeds is
often a very costly process.

In this paper, we introduce a new strategy to compute 
very accurate approximations of feasible points. 
This strategy takes advantage of the Newton 
method for under-constrained systems of equations and inequalities. More precisely, this procedure 
exploits the optimal solution of a linear relaxation of the problem to compute 
efficiently a promising upper bound. First experiments on the Coconuts
benchmarks demonstrate that the combination of this procedure with a safe Branch and Bound algorithm
drastically improves the performances.


\subsection*{The Branch and Bound schema}\label{sec-BB}

The algorithm (see Algorithm \ref{qado-algo}) we describe here is derived from 
the well known Branch and Bound schema introduced by Horst and Tuy
for finding  a global minimizer. Interval analysis techniques are used to  ensure rigorous and safe computations 
whereas constraint programming techniques are used to improve the reduction of the feasible space.

Algorithm \ref{qado-algo} computes enclosers for minimizers and 
safe bounds of the global minimum value within an initial box $\vi{x}$.
Algorithm \ref{qado-algo} maintains two lists : a list \calx{L} of boxes to be processed 
and a list \calx{S} of proven feasible boxes. 
It provides a rigorous encloser $[L, U]$ of the global optimum 
with respect to a tolerance $\epsilon$.

\begin{algorithm}[t]
{\small
\begin{minipage}{.93\linewidth}
{\bf Function} {\tt BB}({\tt IN} $\vi{x}$, $\epsilon$; {\tt OUT} \calx{S}, $[L, U]$)
\begin  {tabbing}
iiii \=   While \= ifif \= ififthen \=  thenthen  \= then \= then  \= \kill
  \% \calx{S}: set of proven feasible points\\
  \% ${\mathbf f}_{\mathbf x}$ denotes the set of possible values for $f$ in ${\mathbf x}$\\
  \% $nbStarts$: number of starting points in the first upper-bounding\\
  $\calx{L} \Xaffect \{{\mathbf x}\}$;\; \;
  $(L, U) \Xaffect (-\infty, +\infty)$;\; \;
  $\calx{S} \Xaffect UpperBounding({\mathbf x}', nbStarts)$;\\
  \Xwhile  $w([L, U]) > \epsilon$ \Xdo \\
  \> ${\mathbf x}' \Xaffect {\mathbf x}''$ such that
         $\lb{{\mathbf f}}_{{\mathbf x}''} = min\{\lb{{\mathbf
             f}}_{{\mathbf x}''}: {\mathbf x}'' \in \calx{L}\}$;\;\; 
     $\calx{L} \Xaffect \calx{L}\backslash{\mathbf x}'$;\;\;
     $\ub{{\mathbf f}}_{{\mathbf x}'} \Xaffect min(\ub{{\mathbf f}}_{{\mathbf x}'}, U)$;\\
  \> ${\mathbf x}' \Xaffect Prune({\mathbf x}')$;\; \;
     $\lb{{\mathbf f}}_{{\mathbf x}'} \Xaffect LowerBound({\mathbf x}')$;\;\; 
     $\calx{S} \Xaffect \calx{S}\cup UpperBounding({\mathbf x}', 1)$;\\ 
  \> \Xif ${\mathbf x}' \neq \emptyset$ \> \> \Xthen \; $({\mathbf x}'_1,
  {\mathbf x}'_2) \Xaffect Split({\mathbf x}')$;\; $\calx{L} \Xaffect
  \calx{L}\cup\{{\mathbf x}'_1, {\mathbf x}'_2\}$;   \\
  \> \Xif $\calx{L}=\emptyset$ \> \> \Xthen  $(L, U) \Xaffect (+\infty, -\infty)$\\
  \> \> \>\Xelse  $(L, U) \Xaffect (min\{\lb{{\mathbf f}}_{{\mathbf
      x}''}: {\mathbf x}'' \in \calx{L}\}, 
                              min\{\ub{{\mathbf f}}_{{\mathbf x}''}:
                              {\mathbf x}'' \in \calx{S}\})$\\ 
  \Xendwhile
\end{tabbing}
\end{minipage}}
\caption{Branch and Bound algorithm}\label{qado-algo}
\end{algorithm}

Algorithm \ref{qado-algo} starts with $UpperBounding({\mathbf x}, nbStarts)$ which computes a set of feasible boxes 
by calling a local search with $nbStarts$ starting points and a proof procedure. 

The box around the local solution is added to $\calx{S}$ if it is proved to contain a feasible point.
At this stage, if the box ${\mathbf x}'$ is empty then,
either it does not contain any feasible point or
its lower bound $\lb{{\mathbf f}}_{{\mathbf x}'}$ is greater than the
current upper bound $U$.  
If ${\mathbf x}'$ is not empty, the box is split along one of the 
variables\footnote{Various heuristics are used to select the variable
  the domain of which has to be split.} of the problem. 

In the main loop, algorithm \ref{qado-algo} selects the box with the lowest lower bound of the
objective function.  
The $Prune$ function applies filtering techniques to reduce the size
of the box ${\mathbf x}'$. In the framework we have implemented,
$Prune$ just uses a 2B-filtering algorithm \cite{Lhomme-1993}. 
Then, $LowerBound({\mathbf x}')$ computes a rigorous lower bound $\lb{{\mathbf f}}_{{\mathbf x}'}$ 
using a linear programming relaxation of the initial problem.
Actually, function $LowerBound$ is based on the linearization
techniques of the {\tt Quad}-framework \cite{Lebbah-Michel-Rueher-2004}. 
$LowerBound$ computes a safe minimizer  $\lb{{\mathbf f}}_{{\mathbf x}'}$ thanks to the techniques 
introduced by Neumaier et al.

Algorithm \ref{qado-algo} maintains the lowest lower bound $L$  
of the remaining boxes $\calx{L}$ and the lowest upper bound $U$ of
proven feasible boxes. 
The algorithm terminates when the space between $U$ and $L$ becomes smaller 
than the given tolerance $\epsilon$.

The Upper-bounding step (see Algorithm \ref{ubfromlb-algo}) performs a multistart strategy 
where a set of $nbStarts$ starting points are provided to a local optimization solver. 
The solutions computed by the local solver are then given to
a function $InflateAndProve$ 
which uses an existence proof procedure
based on the Borsuk test.
However, the proof procedure  may fail to prove the existence of a feasible point within box $\vi{x}_p$.
The most common source of failure is that the ``guess'' provided by the local search
lies too far from the feasible region.


\subsection*{A new upper bounding strategy}\label{ubcorr}

The upper bounding procedure described in the previous section relies on a local
search to provide a ``guessed'' feasible point lying in the neighborhood of a local optima.
However,  the
effects of floating point computation on the provided local optima are hard to predict.
As a result, the local optima might lie outside the feasible region and the
proof procedure might fail to build a proven box around this point.

We propose here a new upper bounding strategy which attempts to take advantage of the 
solution of a linear outer approximation of the problem.
The lower bound process uses such an approximation to compute a safe lower bound
of $\Pb$. When the LP is solved, a solution $x_{LP}$ is always computed and, thus,
available for free. 
This solution being an optimal solution of an outer approximation of $\Pb$,
it lies outside the feasible region. Thus, $x_{LP}$ is not a feasible point.
Nevertheless, $x_{LP}$ may be a good starting point to consider for the following reasons:
\begin{itemize}
\item At each iteration, the branch and bound process splits the domain of the variables.
The smaller the box is, the nearest $x_{LP}$ is from the actual optima of $\Pb$.
\item The proof process inflates a box around the initial guess. This process
may compensate the effect of the distance of $x_{LP}$ from the feasible region.
\end{itemize}
However, while $x_{LP}$ converges to a feasible point, the process might be quite slow.
To speed up the upper bounding process, we have introduced a light weight, though efficient,
procedure which compute a feasible point from a point lying in the neighborhood of the
feasible region. 
This procedure which is called $FeasibilityCorrection$ will be detailed
in the next subsection.

\begin{algorithm}[t]
{\small
\begin{minipage}{.93\linewidth}
{\bf Function} {\tt UpperBounding}({\tt IN} $\vi{x}$, $x^*_{LP}$,
$nbStarts$; {\tt OUT} $\calx{S}'$) 
\begin{tabbing}
iiii \=   While \= ifif \= ififthen \=  thenthen  \= then \= then  \= \kill
  \% $\calx{S}'$: list of proven feasible boxes;\; \; $nbStarts$:
  number of starting points\\ 
  \% $x^*_{LP}$: the optimal solution of the LP relaxation of $\calx{P}(\vi{x})$\\
  $\calx{S}'$ \Xaffect\ $\emptyset$;\;\;
  $x^*_{corr}$ \Xaffect FeasibilityCorrection($x^*_{LP}$); \,\,
  $\vi{x}_p$ \Xaffect InflateAndProve($x^*_{corr}$, $\vi{x}$);\\
  \Xif $\vi{x}_p \neq \emptyset$ \Xthen  $\calx{S}'$ \Xaffect $\calx{S}'\cup\vi{x}_p$\\
  \Xreturn $\calx{S}'$\\
\end{tabbing}
\end{minipage}}
\caption{Upper bounding build from the LP optimal solution $x^*_{LP}$}\label{ubfromlb-algo}
\end{algorithm}

Algorithm \ref{ubfromlb-algo} describes how an upper bound may be build from the
solution of the linear problem used in the lower bounding procedure.


\subsection*{Computing pseudo-feasible points}\label{corr}

This section introduces an adaptation of the Newton method to
under-constrained systems of equations and inequalities  
which provides very accurate approximations of feasible points at a low computational cost.
When the system of equations $\vg(\vx)=\vnull$ is under-constrained  there is a manifold of solutions. 
 $\vl(\vx)$, the linear approximation  is still valid in this situation, 
but the linear system of equations $\vl(\vx)=\vnull$ is now under-constrained, 
and has therefore an affine space of solutions. So we have to choose a solution $\vx^{(1)}$ of the 
linearized equation $\vl(\vx)=\vnull$ among the affine space of solutions. 
As $\vx^{(0)}$ is supposed to be an approximate solution of $\vg(\vx)=\vnull$,
 the best choice is certainly 
the solution of $\vl(\vx)=\vnull$ which is the closest to $\vx^{(0)}$.
This solution can easily be computed with the Moore-Penrose inverse:
$	\vx^{(1)}=\vx^{(0)}-A_{\vg}^+(\vx^{(0)})\vg(\vx^{(0)})$,
where $A_{\vg}^+\in\R^{n\times m}$ is the Moore-Penrose inverse of $A_{\vg}\in\R^{m\times n}$, 
the solution of the linearized equation which minimizes $||\vx^{(1)}-\vx^{(0)}||$. 
Applying previous relation recursively leads to a sequence of vectors 
which converges to a solution close to the initial approximation, 
provided that this latter is accurate enough.

The Moore-Penrose inverse can be computed in several ways: a singular
value decomposition can be used, or in the case where $A_{\vg}$ has
full row rank (which is the case for $A_{\vg}(\vx^{(0)})$   if
$\vx^{(0)}$ is non-singular) the Moore-Penrose inverse can be computed
using $A_{\vg}^+=A_{\vg}^T(A_{\vg}A_{\vg}^T)^{-1}$. 

Inequality constraints are changed to equalities by introducing slack variables:
$	h_j(\vx)\leq0 \iff h_j(\vx)=-s_i^2.$
So the Newton method for under-constrained systems of equations can be applied.

\subsection*{Experiments}\label{expes}

In this Section, we  comment the results of the experiments with our new upper bounding strategy
on a significant set of benchmarks. Detailled results can be found in
the resarch report  ISRN I3S/RR-2008-11-FR\footnote{\small
  http://www.i3s.unice.fr/\%7Emh/RR/2008/RR-08.11-A.GOLDSZTEJN.pdf}).
All 
the benchmarks come from the collection of benchmarks of the Coconuts
project\footnote{See 
http://www.mat.univie.ac.at/{$\tilde{\;}$}neum/glopt/coconut/Benchmark/Benchmark.html.}.
We have selected $35$  benchmarks
 where Icos succeeds to find the global
minimum while relying on an unsafe local search.
We did compare our new upper bounding strategy with  the following upper
bounding strategies:
\begin{itemize}
\item[S1:] This strategy directly uses the guess from the local search, i.e.
     this strategy uses a simplified version of algorithm \ref{qado-algo} where
     the proof procedure has been dropped.
     As a consequence, it does not suffer from the difficulty to prove the existence of a feasible point.
     However, this strategy is unsafe and may produce wrong results.
\item[S2:] This strategy attempts to prove the existence of a feasible point within a box
     build around the local search guess. 
     Here, all provided solutions are safe and the solving process is
     rigorous.
\item[S3:] Our upper bounding strategy where the upper
     bounding relies on the optimal solution of the problem linear
     relaxation to  
     build a box proved to hold a feasible point.
     A call to the correction procedure attempts to
     compensate the effect of the outer approximation.
\item[S4:] This strategy  applies the correction procedure to the
  output of the local search (to improve the approximate solution given by a
     local search).
\item[S5:] This strategy mainly differs from S3 by the fact that it does not
     call the correction procedure
\end{itemize}
S3, our new upper bounding strategy is the best strategy:
31 benchmarks are now solved within the 30s time out; moreover, almost
all benchmarks are solved in much less time and with a greater amount
of proven solutions. 
This new strategy improves drastically the performance of the upper bounding procedure and
 competes well with a local search.  

Current work aims at improving and generalizing this framework and its implementation. 

\bibliographystyle{plain}
\bibliography{cp08_short}

\end{document}